\documentclass[12pt]{article}

\newcommand{\frak}{\cal}
\newcommand{\tree}{\mathop{\rm Tree}}
\newcommand{\forest}{\mathop{\rm Forest}}
\newcommand{\delroot}{\mathop{\rm DelRoot}}
\newcommand{\hangon}[1]{\mathbin{{\leftharpoonup}_{\bf #1}}}

\newcommand{\nodes}{\mathop{\rm Nodes}}
\newcommand{\troot}{\mathop{\rm root}}
\newcommand{\card}{\mathop{\rm Card}}
\newcommand {\frX}{{\frak X}}
\newcommand {\frT}{{\frak T}}
\newcommand {\frOT}{{\frak O\frak T}}
\newcommand {\frHOT}{{\frak H\frak O\frak T}}
\newcommand {\frLT}{{\frak L\frak T}}
\newcommand {\frLOT}{{\frak L\frak O\frak T}}
\newcommand {\frLHOT}{{\frak L\frak H\frak O\frak T}}
\newcommand {\frGT}{{\frak G\frak T}}
\newcommand {\frLGT}{{\frak L\frak G\frak T}}
\newcommand{\Lmap}{{\rm L}}
\newcommand{\Rmap}{{\rm R}}
\newcommand{\Mmap}{{\rm M}}
\newcommand{\Dmap}{{\rm D}}
\newcommand{\bfd}{{\bf d}}
\newcommand{\bfe}{{\bf e}}
\newcommand{\bff}{{\bf f}}
\newcommand{\bfg}{{\bf g}}
\newcommand{\kX}{{k\{\frX\}}}
\newcommand{\kT}{{k\{\frT\}}}
\newcommand{\kFLOT}{k{<}\frLOT_1{>}}
\newcommand{\kIX}{k{<}I(\frX){>}}
\newcommand{\kLOT}{k\{\frLOT\}}
\newcommand{\kHOT}{k\{\frHOT\}}
\newtheorem{axiom}{Axiom}
\newtheorem{lemma}{Lemma}[section]

\newtheorem{thm}[lemma]{Theorem}
\newtheorem{corollary}[lemma]{Corollary}
\newtheorem{example}[lemma]{Example}
\newcommand{\pf}{\medbreak\noindent{\sc Proof:}\ }

\author{Robert Grossman\thanks{The
first author is a National Science Foundation Postdoctoral
Research Fellow.} \\ University of California, Berkeley
\and Richard G. Larson\thanks{This
paper was written while the second
author was on sabbatical leave at the University of California,
Berkeley.}\\University of Illinois at Chicago }
\title{Hopf-algebraic structures\\
of families of trees}

\date{February, 1987}

\begin{document}
\maketitle

\medbreak
{\noindent\bf
This is a draft of a paper which subsequently appeared
in the Journal of Algebra, Volume 26, 1989, pp. 184--210.}
\medbreak

\section{Introduction}

In this paper we describe Hopf algebras which are associated with
certain families of trees.
These Hopf algebras originally arose in a natural fashion:
one of the authors~\cite{Grossman} was investigating data structures
based on trees, which could be used to efficiently compute certain
differential operators.
Given data structures such as trees which can be multiplied, and which
act as
higher-order derivations on an algebra, one expects to find a Hopf
algebra of some sort.
We were pleased to find that not only was there a Hopf algebra
associated with these data structures,
but that
it could be used to give new proofs of enumerations of such
objects as rooted trees and ordered rooted trees.
Previous work applying Hopf algebras to combinatorial objects
(such as~\cite{JoniRota}, \cite{NiSw} or~\cite{Bodo}) has concerned
itself
with algebraic structures on polynomial algebras and on partially
ordered sets, rather than on trees themselves.
We hope that these constructions will also provide insight for
the algebra of data structures.

The Hopf algebras which we construct are all cocommutative graded
connected Hopf algebras.
This allows us to apply the Milnor-Moore Theorem (Theorem~\ref{MMThm})
and the Poincar\'e-Birkhoff-Witt Theorem (Theorem~\ref{PBWThm}) to get
precise information on the structure of these Hopf algebras.
We illustrate our construction, and its application, by
sketching a proof of Cayley's enumeration~\cite{Cayley} of
finite rooted trees.

We now describe how to
construct the Hopf algebra $\kT$ which has as basis all finite
rooted trees.
The grading on $\kT$ is given as follows:
if the tree $t$ has $n+1$ nodes, then $t\in\kT_n$.
If $t_1$ and $t_2$ are trees, the product $t_1\cdot t_2$ is the sum of
the trees formed by attaching the children of the root of $t_1$ to the
nodes of $t_2$ in all possible ways.
If $t$ is a tree, the coproduct $\Delta(t)$ is the sum of all terms
$t_1\otimes t_2$, where the children of the root of $t_1$ and the
children of the root of $t_2$ range over all possible partitions of the
children of the root of $t$.
This definition of the coproduct is very similar to the
definition of the coproduct in the {\em placement coalgebra\/} described
in~\cite{JoniRota}.
Both here and in~\cite{JoniRota},
the coproduct of a structure is the sum of all
terms which are
the tensor product of the two pieces resulting from decomposing
in all possible ways
the structure into two disjoint substructures.

We show
in Sections~\ref{CombSect} and~\ref{HopfAlgSect} that this
product and coproduct gives a cocommutative graded connected Hopf
algebra on the vector space which has as basis the isomorphism classes
of finite rooted trees.
Specifically, in Section~\ref{CombSect} we give axioms which a family of
trees must satisfy to support a Hopf algebra structure.
In Section~\ref{HopfAlgSect} we define the Hopf algebra associated with
such a family.
In our example, it is straightforward to show (Theorem~\ref{prims}) that
$P(\kT)$,
the space of primitive elements in $\kT$, has as basis the set of trees
whose root has exactly one child.
Therefore we have a one-one correspondence between a basis for
$\kT_{n-1}$ and a basis for $P(\kT)_n$:
any tree with $n$ nodes corresponds to the tree with $n+1$ nodes formed
by creating a new root and linking the root of the original tree to it
as a child.
If we let $t_n$ be the number of trees with $n$ nodes,
$a_n=\dim \kT_n$, and $p_n=\dim P(\kT)_n$, we have that
$t_{n+1}=a_n$ and $p_n=a_{n-1}$.
The Milnor-Moore Theorem and the Poincar\'e-Birkhoff-Witt Theorem imply
that
\begin{eqnarray*}
a_{n} & = & \sum_{m_1 + 2m_2 +\cdots+rm_r = n}
{m_1 + p_1 - 1\choose m_1}\cdots{m_r + p_r - 1\choose m_r}\\
      & = & \sum_{m_1 + 2m_2 +\cdots+rm_r = n}
{m_1 + a_0 - 1\choose m_1}\cdots{m_r + a_{r-1} - 1\choose m_r},
\end{eqnarray*}
which implies that
\[
t_{n+1} = \sum_{m_1 + 2m_2 +\cdots+rm_r = n}
{m_1 + t_1 - 1\choose m_1}\cdots{m_r + t_r - 1\choose m_r},
\]
which implies, using
\[
{-r\choose k} = {(-1)}^k{r+k-1\choose k},
\]
Cayley's result that if
$T(z) = \sum_{n = 1}^\infty t_n z^n$
then
\[
T(z) = z \prod_{n=1}^\infty {(1 - z^n)}^{-t_n}.
\]

In Section~\ref{LTSect} we give the proof we have just sketched, in full
detail, for labeled (or colored) trees.
In Section~\ref{LOTSect} we discuss the algebraic structure of the Hopf
algebras constructed for the family of labeled ordered trees (that is,
for labeled or colored trees in which the children of each root
are linearly ordered).
In Section~\ref{XSect} we discuss the algebraic structure of the Hopf
algebras constructed for the family of heap-ordered trees, and for
related families.
Heap-ordered trees (see~\cite{Tarjan} for details) are used as data
structures in computer science and symbolic algebra.
In this section, we also give a brief discussion of the category of
families of trees.

Throughout this paper, the field $k$ will have characteristic~0.


\section{The structure of families of rooted trees}\label{CombSect}

In this section we describe a structure we impose on families of
rooted trees.
By a tree we mean a nonempty finite rooted tree, and by
a forest we mean a finite family of finite rooted trees, possibly
empty.
The mathematical structure $\frX$ consists of a family
$\tree(\frX)$ of trees, and a family $\forest(\frX)$  of
forests, both possibly with additional structure,
together with four operations.
In most of the examples we
consider, the additional structure consists of
orderings or labels.
As we proceed, we impose axioms requiring that the
operations preserve this additional structure.
The operations are:
\begin{itemize}
\item an operation $\delroot$ which maps $\tree(\frX)$ to
$\forest(\frX)$.
This sends a tree $t \in \tree(\frX)$ to the forest obtained when
the root of $t$ is deleted.  Note that the exact definition of
$\delroot$
depends  upon the additional structure we impose on $\tree(\frX)$ and
$\forest(\frX)$.
\item an operation $\nodes$ which maps $\tree(\frX)$ to {\bf Sets}.
This sends a tree $t \in \tree(\frX)$ to the set consisting of the
nodes of the underlying tree of $t$.
\item a restriction operation which maps subforests of $\forest(\frX)$
to $\forest(\frX)$.
If $V \in \forest(\frX)$, and if $U \subseteq V$ is a
subforest,
then we must specify a way to impose the additional structure of
$\frX$ on the forest $U$ in a manner reflecting the structure of $V$.
We denote the resulting element of $\forest(\frX)$ by
$V|U$.
\item an attachment operation which is a map from
$\forest(\frX) \times \tree(\frX)$ to $\tree(\frX)$.
We denote this operation by $\leftharpoonup$.
If $t \in \tree(\frX)$,
$U \in \forest(\frX)$, and $\bfd : U \rightarrow \nodes(t)$,
then $t\hangon{d}U \in \tree(\frX)$.
Intuitively, what $\leftharpoonup$ does is
the following: it forms a new tree by linking the root of each tree $u$
in the forest $U$ to the node $\bfd(u)$ of $t$, in a manner which
preserves the additional structure of $\frX$.
\end{itemize}

We give some examples of specific families of trees we consider in this
paper.
For each of these examples, it is easy to check that the nine axioms
below are satisfied.

\begin{example}
The simplest example is $\frT$, the family of trees without any
additional structure.
The set $\tree(\frT)$ is the set of finite rooted trees.
The set $\forest(\frT)$ is the set of finite forests of finite rooted
trees.
The map $\delroot$ sends each tree into the forest formed by deleting
its root.
If $V\in\forest(\frT)$ and $U\subset V$ is a subforest, then $V|U$ is
$U$.
If $t$ is a finite rooted tree, $U$ is a finite forest of finite rooted
trees, and $\bfd : U\rightarrow\nodes(t)$ is a function,
then $t\hangon{d}U$ is the tree formed by linking the root of each tree
in $U$ to the node $\bfd(u)$ of $t$.
\end{example}

Most of the axioms we give for the family $\frX$ will
consist
of assertions that facts we explicitly prove for $\frT$ hold for
$\frX$.

\begin{example}
The family $\frOT$ of ordered trees.
The set $\tree(\frOT)$ is the set of finite rooted trees, together with
total orderings on the sets of children of each node.
The set $\forest(\frOT)$ is the set of finite forests of finite rooted
trees, together with total orderings on the sets of children of each
node, and a total ordering on each forest.
The map $\delroot$ sends each tree into the forest formed by deleting
its root.
The ordering on the forest is the one arising from the ordering on the
children of the root of the original tree.
If $V\in\forest(\frOT)$ and $U\subset V$ is a subforest, then $V|U$ is
$U$, with the same orderings on the sets of children of each node as in
$V$, and with the ordering on the trees in $U$ induced as a subset of
the trees of $V$.
If $t$ is an ordered finite rooted tree, $U$ is an ordered finite forest
of ordered
finite rooted trees, and $\bfd : U\rightarrow\nodes(t)$ is a function,
then $t\hangon{d}U$ is the tree formed by linking the root of each tree
in $U$ to the node $\bfd(u)$ of $t$.
The ordering on the children of a node $\alpha$ of
$t\hangon{d}U$ is given as follows:
the children newly linked to $\alpha$ follow all of the
original children of $\alpha$ in
the new order, and keep the order among themselves which was induced by
the ordering of the forest $U$;
the original children of $\alpha$ keep their original order.
\end{example}

Note that in the above example, nodes which do not have a common parent
are not related by any order.
Also note that we have orderings on subsets of sets of nodes of
trees, and also have orderings on forests of trees.

\begin{example}
The family $\frHOT$ of heap-ordered trees.
The set $\tree(\frHOT)$ is the set of finite rooted trees, together with
a total ordering on the set of all nodes of each tree,
such that each node precedes all of its children in the
ordering.
The set $\forest(\frHOT)$ is the set of finite forests of finite rooted
trees, together with a total ordering on the set of all the nodes of
each forest, such that each node precedes all of its children in the
ordering.
The map $\delroot$ sends each tree into the forest formed by deleting
its root.
The ordering on the nodes of the forest is the one arising from the
ordering on the nodes of the original tree.
If $V\in\forest(\frHOT)$ and $U\subset V$ is a subforest, then $V|U$ is
$U$, with the ordering on the nodes of $U$ being induced as a subset of
the nodes of $V$.
If $t\in\tree(\frHOT)$, $U\in\forest(\frHOT)$,
and $\bfd : U\rightarrow\nodes(t)$ is a function,
then $t\hangon{d}U$ is the tree formed by linking the root of each tree
in $U$ to the node $\bfd(u)$ of $t$.
The new ordering on the nodes of
$t\hangon{d}U$ is given as follows:
the nodes in the forest $U$ follow all of the
nodes of the tree $t$ in
the new order, and keep the order among themselves which was induced by
the ordering of the nodes of the trees of the forest $U$;
the original nodes in $t$ keep their original order.
\end{example}

\begin{example}\label{LTEx}
The family $\frLT(E_1$, \ldots, $E_M)$ of labeled trees, where
$\{E_1$, \ldots, $E_M\}$ is a set of formal symbols.
The set $\tree(\frLT(E_1$, \ldots, $E_M))$ is the set of finite rooted
trees,
together with an assignment of an element of $\{E_1$, \ldots, $E_M\}$
to each node, other than the root, of each tree.
The set $\forest(\frLT(E_1$, \ldots, $E_M))$
is the set of finite forests of finite rooted
trees,
together with an assignment of an element of $\{E_1$, \ldots, $E_M\}$
to each node in the forest.
The map $\delroot$ sends each tree into the forest formed by deleting
its root.
The assignment of symbols to the nodes is not changed.
If $V\in\forest(\frLT(E_1$, \ldots, $E_M))$ and $U\subset V$
is a subforest, then $V|U$ is $U$ with the same assignment
of formal symbols.
If $t\in\tree(\frLT)$, $U\in\forest(\frLT)$,
and $\bfd : U\rightarrow\nodes(t)$ is a function,
then $t\hangon{d}U$ is the tree formed by linking the root of each tree
in $U$ to the node $\bfd(u)$ of $t$.
The assignment of symbols to the nodes is not changed.

The families $\frLOT$ of labeled ordered trees,
and $\frLHOT$ of labeled heap-ordered trees, are defined analogously.
\end{example}

Labeled trees are called {\em colored} trees by some authors.
Note that for $\frX=\frLT$, $\frLOT$, or $\frLHOT$, the trees in
$\tree(\frX)$ do not have their roots labeled, but the trees in forests
in $\forest(\frX)$ have all of their nodes, including their roots,
labeled.
The above examples all occur in applications
to analysis and to data structures used in
symbolic algebra.

\medbreak
The first axiom avoids degenerate families $\frX$.

\begin{axiom}
Every finite rooted tree occurs as the underlying tree of some element
of $\tree(\frX)$.
Every finite forest of finite rooted trees appears as the underlying
forest of some element of $\forest(\frX)$.
\end{axiom}

The next three axioms describe how the operation $|$ relates to the
structure of elements of $\forest(\frX)$.

\begin{axiom}\label{TrivRestr}
Let $V\in\forest(\frX)$.
Then $V|V=V$.
\end{axiom}

\begin{axiom}
Let $V\in\forest(\frX)$, and let $U\subset V$ be a subforest.
Then the underlying forest of $V|U$ is $U$.
\end{axiom}

\begin{axiom}\label{CoAssocAxiom}
Let $W\in\forest(\frX)$, and let $U\subset V\subset W$ be subforests.
Then
\[
W|U=(W|V)|U.
\]\end{axiom}

The next axiom describes how the linking operation $\leftharpoonup$
interacts with the structure of the elements of $\tree(\frX)$ and
$\forest(\frX)$.

\begin{axiom}
Let $t\in\tree(\frX)$, $U\in\forest(\frX)$,
and $\bfd : U\rightarrow\nodes(t)$.
Then $t \hangon{d} U \in \tree(\frX)$.
The underlying tree of $t \hangon{d} U$ is the tree formed by linking
the root of each tree $u\in U$ to the node $\bfd(u)$ of $t$.
\end{axiom}

\begin{axiom}\label{UnitAxiom}
There exists a unique $e\in\tree(\frX)$ such that the underlying tree of
$e$ has exactly one node.

There exists a unique $\emptyset\in\forest(\frX)$ such that the
underlying forest of $\emptyset$ is the empty set.
\end{axiom}

For any $U\in\forest(\frX)$ there is a unique function
\[
{\bf 1} : U\rightarrow\nodes(e).
\]

For any $t\in\tree(\frX)$ there is a unique function
\[
{\bf 0} : \emptyset\rightarrow\nodes(t).
\]

\begin{axiom}\label{HangonAxiom}
If $t\in\tree(\frX)$, then
\[
t = e \hangon{1} \delroot(t),
\]
and
\[
t = t \hangon{0} \emptyset.
\]
If $U\in\forest(\frX)$, then
\[
U = \delroot(e \hangon{1} U).
\]
\end{axiom}

The previous axiom implies that there is redundancy in the way we have
defined our structure $\frX$:
$\tree(\frX)$ and $\forest(\frX)$ are essentially isomorphic via the
maps $t\mapsto\delroot(t)$ and $U\mapsto(e\hangon{1}U)$.
However, we have chosen to retain both $\tree(\frX)$ and
$\forest(\frX)$ for purposes of clarity and intuitiveness.

If we take three finite rooted trees $t_1$, $t_2$, and $t_3$, there are
two ways to form a tree:
with the children of the root of $t_1$ linked to the nodes
of $t_2$,
and with the children of the root of $t_2$ linked to the nodes of $t_3$.
We show that these two ways are the same.
(This is essentially an associativity condition.)

Suppose that
\[
\bfd : \delroot(t_1)\rightarrow\nodes(t_2),
\]
and
\[
\bfe : \delroot(t_2\hangon{d}\delroot(t_1))\rightarrow\nodes(t_3)
\]
are linking functions.
Then we can form the tree
\[
t_3\hangon{e}\delroot(t_2\hangon{d}\delroot(t_1)).
\]
From these data we can construct linking functions
\[
\bff : \delroot(t_2)\rightarrow\nodes(t_3)
\]
and
\[
\bfg : \delroot(t_1)\rightarrow\nodes(t_3\hangon{f}\delroot(t_2))
\]
as follows.
Every tree in the forest
$\delroot(t_2)$ corresponds to a tree in the forest
$\delroot(t_2\hangon{d}\delroot(t_1))$.
Let
\refstepcounter{lemma}\label{iota}
$$
\iota : \delroot(t_2)\rightarrow\delroot(t_2\hangon{d}\delroot(t_1))
\eqno{\rm (\thelemma)}
$$
be the injection which gives this correspondence.
The function $\bf f$ is given by
\[
\bff=\bfe\circ \iota.
\]
The function $\bf g$ is given by
\[
\bfg(s) = \left\{
\begin{array}{ll}
\bfd(s) & \mbox{if $\bfd(s)\not=\troot(t_2)$}\\
\bfe(s) & \mbox{if $\bfd(s)=\troot(t_2)$.}
\end{array} \right.
\]
Note that we identify $\delroot(t_2\hangon{d}\delroot(t_1))$
with the disjoint union
\[
\iota(\delroot(t_2)) \cup \{\,u\in\delroot(t_1)\mid\bfd(u)=\troot(t_2)\,\},
\]
and identify $\nodes(t_3\hangon{f}\delroot(t_2))$
with the disjoint union of $\nodes(t_3)$ and
$\nodes(t_2)\backslash\{\troot(t_2)\}$.
($X\backslash Y$ denotes the set-theoretic difference of $X$ and $Y$.)
We denote the map sending the pair $(\bfd,\bfe)$
to the pair $(\bff,\bfg)$ as follows:
\[
\Rmap(\bfd,\bfe)=(\bff,\bfg).
\]
Suppose now that
\[
\bff : \delroot(t_2)\rightarrow\nodes(t_3),
\]
and
\[
\bfg : \delroot(t_1)\rightarrow\nodes(t_3\hangon{f}\delroot(t_2))
\]
are linking functions.
Then we can form the tree
\[
(t_3\hangon{f}\delroot(t_2))\hangon{g}\delroot(t_1).
\]
From these data we can construct linking functions
\[
\bfd : \delroot(t_1)\rightarrow\nodes(t_2)
\]
and
\[
\bfe : \delroot(t_2\hangon{d}\delroot(t_1))\rightarrow\nodes(t_3)
\]
as follows.
The function $\bf d$ is given by
\[
\bfd(s) = \left\{
\begin{array}{ll}
\bfg(s) & \mbox{if $\bfg(s)\in\nodes(t_2)\backslash\{\troot(t_2)\}$}\\
\troot(t_2) & \mbox{otherwise.}
\end{array} \right.
\]
The function $\bf e$ is given by
\[
\bfe(s) = \left\{
\begin{array}{ll}
\bff(s') & \mbox{if $s=\iota(s')$ with $s'\in\delroot(t_2)$}\\
\bfg(s)  & \mbox{if $s\in\delroot(t_1)$
   and $\bfg(s)\in\nodes(t_3)$.}
\end{array} \right.
\]
Note that we identify $\delroot(t_2\hangon{d}\delroot(t_1))$
with the disjoint union
\[
\iota(\delroot(t_2)) \cup \{\,u\in\delroot(t_1)\mid\bfd(u)=\troot(t_2)\,\},
\]
and identify
$\nodes(t_3\hangon{f}\delroot(t_2))$
with the disjoint union of $\nodes(t_3)$ and
$\nodes(t_2)\backslash\{\troot(t_2)\}$.
We denote the map sending the pair $(\bff,\bfg)$
to the pair $(\bfd,\bfe)$ as follows:
\[
\Lmap(\bff,\bfg)=(\bfd,\bfe).
\]

The following lemma asserts that $\Lmap$ and $\Rmap$ are inverses
of each other.

\begin{lemma}\label{AssocPairsLemma}
Let $t_1$, $t_2$, $t_3\in\tree(\frX)$, and let
\begin{eqnarray*}
\bfd & : & \delroot(t_1)\rightarrow\nodes(t_2)\\
\bfe & : & \delroot(t_2\hangon{d}\delroot(t_1))\rightarrow\nodes(t_3)
\end{eqnarray*}
be functions.
Then
\[
\Lmap\circ\Rmap(\bfd,\bfe) = (\bfd,\bfe).
\]
Let
\begin{eqnarray*}
\bff & : & \delroot(t_2)\rightarrow\nodes(t_3)\\
\bfg & : & \delroot(t_1)\rightarrow\nodes(t_3\hangon{f}\delroot(t_2))
\end{eqnarray*}
be functions.
Then
\[
\Rmap\circ\Lmap(\bff,\bfg) = (\bff,\bfg).
\]
\end{lemma}
\pf We first show that
$\Lmap\circ\Rmap(\bfd,\bfe) = (\bfd,\bfe)$.
Let $\Rmap(\bfd,\bfe) = (\bff,\bfg)$, and
let $\Lmap(\bff,\bfg) = (\bfd',\bfe')$.
To show that $\bfd' = \bfd$ we consider two cases:
if $\bfd(s)\not=\troot(t_2)$, then $\bfg(s)=\bfd(s)$,
and $\bfd'(s)=\bfg(s)$;
if $\bfd(s)=\troot(t_2)$, then $\bfg(s)=\bfe(s)\in
\nodes(t_3)$, and $\bfd'(s)=\troot(t_2)$.
Therefore, in either case, $\bfd(s)=\bfd'(s)$.

We next show that $\bfe'=\bfe$.
Suppose $s\in\delroot(t_2\hangon{d}\delroot(t_1))$.
There are two cases to consider.
In the first case, $s=\iota(s')$ for $s'\in\delroot(t_2)$.
Then $\bfe'(s)=\bff(s')=\bfe\circ \iota(s')=\bfe(s)$.
In the second case, $s\in\delroot(t_1)$, and $\bfd(s)=\troot(t_2)$.
Then $\bfg(s)=\bfe(s)\in\nodes(t_3)$, so
$\bfe'(s)=\bfg(s)=\bfe(s)$.
Therefore, in either case,  $\bfe'(s)=\bfe(s)$.

We finally show that
$\Rmap\circ\Lmap(\bff,\bfg) = (\bff,\bfg)$.
Let $\Lmap(\bff,\bfg) = (\bfd,\bfe)$, and
let $\Rmap(\bfd,\bfe) = (\bff',\bfg')$.
If $s\in\delroot(t_2)$, then $\bff'(s)=\bfe\circ \iota(s)=\bff(s)$.
We show that $\bfg'=\bfg$ by considering two cases.
In the first case, $\bfg(s)\in\nodes(t_2)\backslash\{\troot(t_2)\}$.
Then $\bfd(s)=\bfg(s)$, and $\bfd(s)\not=\troot(t_2)$.
Therefore $\bfg'(s)=\bfd(s)=\bfg(s)$.
In the second case, $\bfg(s)\in\nodes(t_3)$.
Then $\bfd(s)=\troot(t_2)$.
Now, since $s\in\delroot(t_1)$ and $\bfe(s)\in\nodes(t_3)$, we have that
$\bfg'(s)=\bfe(s)=\bfg(s)$.
Therefore $\bfg'=\bfg$.
This completes the proof of the lemma.
\medbreak

\begin{lemma}\label{EqualTreeLemma}
Let $t_1$, $t_2$, $t_3\in\tree(\frX)$, let
\begin{eqnarray*}
\bfd & : & \delroot(t_1)\rightarrow\nodes(t_2)\\
\bfe & : & \delroot(t_2\hangon{d}\delroot(t_1))\rightarrow\nodes(t_3)
\end{eqnarray*}
be functions, and let
\[
(\bff,\bfg)=\Rmap(\bfd,\bfe).
\]
Then the underlying trees of
\[
t_3\hangon{e}\delroot(t_2\hangon{d}\delroot(t_1))
\]
and
\[
(t_3\hangon{f}\delroot(t_2))\hangon{g}\delroot(t_1)
\]
are equal.
\end{lemma}

\pf Let
\[
u = t_3\hangon{e}\delroot(t_2\hangon{d}\delroot(t_1))
\]
and
\[
v = (t_3\hangon{f}\delroot(t_2))\hangon{g}\delroot(t_1).
\]
The sets of nodes of the underlying trees of $u$ and $v$ are the same:
the disjoint union of $\nodes(t_3)$,
$\nodes(t_2)\backslash\{\troot(t_2)\}$,
and
$\nodes(t_1)\backslash\{\troot(t_1)\}$.
Let $\alpha$ be one of these nodes.
We consider the three possible cases.
If $\alpha\in\nodes(t_3)$,
then the node $\alpha$ has the same
parent in both $u$ and $v$ as it has in $t_3$.
If $\alpha\in\nodes(t_2)\backslash\{\troot(t_2)\}$,
then there are two subcases to consider.
If $\alpha$ is not a child of the root of $t_2$,
then the node $\alpha$ has the same
parent in both $u$ and $v$ as it has in $t_2$.
If $\alpha$ is a child of the root of $t_2$, then $\alpha$ is the root of some
tree $s\in\delroot(t_2)$.
In $u$, the parent of $\alpha$ is $\bfe\circ\iota(s)$.
In $v$, the parent of $\alpha$ is $\bff(s)=\bfe\circ\iota(s)$.
Therefore, in both subcases, $\alpha$ has the same node as parent in both
$u$ and $v$.
If $\alpha\in\nodes(t_1)\backslash\{\troot(t_1)\}$,
then there are two subcases to consider.
If $\alpha$ is not a child of the root of $t_1$,
then the node $\alpha$ has the same
parent in both $u$ and $v$ as it has in $t_1$.
If $\alpha$ is a child of the root of $t_1$, then $\alpha$ is the root of some
tree $s\in\delroot(t_1)$.
In $u$, if $\bfd(s)\not=\troot(t_2)$, then the parent of $\alpha$ is
$\bfd(s)$;
if $\bfd(s)=\troot(t_2)$, then the parent of $\alpha$ is $\bfe(s)$.
In $v$, the parent of $\alpha$ is $\bfg(s)$.
But $\bfg(s)=\bfd(s)$ if $\bfd(s)\not=\troot(t_2)$,
and $\bfg(s)=\bfe(s)$ if $\bfd(s)=\troot(t_2)$.
Therefore, in both subcases, $\alpha$ has the same node as parent in both
$u$ and $v$.
So in all cases, the node $\alpha$ has the same parent in the underlying
trees of both $u$ and $v$.
It follows that both $u$ and $v$ have the same underlying tree.
This completes the proof of the lemma.
\medbreak

The next axiom says that the equality of trees in
Lemma~\ref{EqualTreeLemma} holds in the family $\frX$.
\begin{axiom}\label{EqualTreeAxiom}
Let $\frX$ be a family of trees, let $t_1$, $t_2$, $t_3\in\tree(\frX)$,
let
\begin{eqnarray*}
\bfd & : & \delroot(t_1)\rightarrow\nodes(t_2)\\
\bfe & : & \delroot(t_2\hangon{d}\delroot(t_1))\rightarrow\nodes(t_3)
\end{eqnarray*}
be functions, and let
\[
(\bff,\bfg)=\Rmap(\bfd,\bfe).
\]
Then
\[
t_3\hangon{e}\delroot(t_2\hangon{d}\delroot(t_1))=
(t_3\hangon{f}\delroot(t_2))\hangon{g}\delroot(t_1).
\]
\end{axiom}

We now develop the material leading up to an axiom necessary for the
``Hopf condition'' on the Hopf algebra we will associate with the
family $\frX$.

Fix $t\in\tree(\frX)$ and  $X\in\forest(\frX)$.
Suppose that
$\bfd : X \rightarrow \nodes(t)$ is a linking function,
and that
$U\subset\delroot(t\hangon{d}X)$ is a subforest.
We define subforests
\[
V\subset\delroot(t)
\]
and
\[
W\subset X
\]
as follows.
Let
\[
V=\{\,v\in\delroot(t) \mid \iota(v)\in U\,\}
\]
and
\[
W=\left\{\,w\in X \,\left|\,
\begin{array}{l}
\mbox{\rm $\bfd(w)=\troot(t)$ and $w\in U$, or}\\
\mbox{$\bfd(w)\in\nodes(v)$ for some $v\in V$}
\end{array}
\,\right\}\right.
\]
and define linking functions
\[
\bfe : W \rightarrow \nodes(e\hangon{1}\delroot(t)|V)
\]
and
\[
\bff : X\backslash W \rightarrow
\nodes(e\hangon{1}\delroot(t)|(\delroot(t)\backslash V))
\]
as follows:
\[
\bfe(w) = \left\{
\begin{array}{ll}
\troot(e) & \mbox{if $\bfd(w)=\troot(t)$}\\
\bfd(w)   & \mbox{if $\bfd(w)\in\nodes(v)$ for some $v\in V$}
\end{array} \right.
\]
and
\[
\bff(w) = \left\{
\begin{array}{ll}
\troot(e) & \mbox{if $\bfd(w)=\troot(t)$}\\
\bfd(w)   & \mbox{if $\bfd(w)\in\nodes(v)$ for some
$v\in\delroot(t)\backslash V$.}
\end{array} \right.
\]
Note that the set over which $U$ ranges depends on $\bfd$,
and the sets over which $\bfe$ and $\bff$ range depend on $V$ and
$W$.

We denote the map sending $(\bfd, U)$ to $(V, W, \bfe,\bff)$
as follows:
\[
\Mmap(\bfd, U)=(V, W, \bfe, \bff).
\]
Suppose now that $V \subset \delroot(t)$ and $W \subset X$
are subforests, and that
\begin{eqnarray*}
\bfe & : & W \rightarrow \nodes(e\hangon{1}\delroot(t)|V)\\
\bff & : & X\backslash W \rightarrow \nodes\bigl(e\hangon{1}
             \delroot(t)|(\delroot(t)\backslash V)\bigr)
\end{eqnarray*}
are linking functions.
Let $\bfd : X \rightarrow \nodes(t)$
be the linking function defined as follows:
\[
\bfd(x) = \left\{
\begin{array}{ll}
\bfe(x)  & \mbox{if $x\in W$}\\
\bff(x) & \mbox{if $x\notin W$,}
\end{array} \right.
\]
and let
\[
U=\iota(V)\cup\{\,x\in W \mid \bfd(x)=\troot(t)\,\}.
\]
(Recall that $\iota$ was defined in Equation~(\ref{iota}).)
We denote the map associating $(\bfd, U)$ to $(V, W, \bfe, \bff)$
as follows:
\[
\Dmap(V, W, \bfe, \bff) = (\bfd, U).
\]

The following lemma asserts that $\Mmap$ and $\Dmap$ are inverses
of each other.

\begin{lemma}\label{HopfCondIndex}
Suppose $t\in\tree(\frX)$ and $X\in\forest(\frX)$.
Let
\[
\bfd : X \rightarrow \nodes(X)
\]
be a function, and let $U\subset\delroot(t\hangon{d}X)$ be a
subforest.
Then
\[
\Dmap\circ\Mmap(\bfd, U) = (\bfd, U).
\]
Let $V\subset\delroot(t)$ and $W\subset X$ be subforests, and
\begin{eqnarray*}
\bfe & : & W \rightarrow \nodes(e\hangon{1}\delroot(t)|V)\\
\bff & : & X\backslash W \rightarrow \nodes\bigl(e\hangon{1}
             \delroot(t)|(\delroot(t)\backslash V)\bigr)
\end{eqnarray*}
be functions.
Then
\[
\Mmap\circ\Dmap(V, W, \bfe, \bff) = (V, W, \bfe, \bff).
\]
\end{lemma}

\pf We first show that $\Dmap\circ\Mmap(\bfd, U) = (\bfd, U)$.
Let $\Mmap(\bfd, U) = (V, W, \bfe, \bff)$, and let
$\Dmap(V, W, \bfe, \bff) = (\bfd', U')$.
We first show that $\bfd' = \bfd$.
Using the identification $t = e\hangon{1}\delroot(t)$ from
Axiom~\ref{HangonAxiom}, we have that $\troot(t)=\troot(e)$ in
$\nodes(t) = \nodes\bigl(e\hangon{1}\delroot(t)\bigr)$.
For $x\in X$, we have that $\bfd'(x)$ is $\bfe(x)$ or $\bff(x)$,
which is $\troot(e)$ if $\bfd(x)=\troot(t)$, and $\bfd(x)$ otherwise.
Therefore $\bfd'(x) = \bfd(x)$.
Since
\[
\delroot(t\hangon{d}X) = \iota(\delroot(t)) \cup
\{\,x\in X \mid \bfd(x) = \troot(t)\,\}
\]
as forests, we have
\begin{eqnarray*}
U' & = & \iota(V) \cup \{\,x\in W \mid \bfd(x) = \troot(t)\,\}\\
   & = & \bigl(\iota(\delroot(t)) \cap U\bigr) \cup
            \bigl(\{\,x\in X \mid \bfd(x) = \troot(t)\,\} \cap U\bigr)\\
   & = & \bigl(\iota(\delroot(t)) \cup \{\,x\in X \mid \bfd(x) =
            \troot(t)\,\}\bigr) \cap U\\
   & = & \delroot(t\hangon{d}X)\cap U\\
   & = & U.
\end{eqnarray*}
Therefore $\Dmap\circ\Mmap(\bfd, U) = (\bfd, U)$.

We now show that
\[
\Mmap\circ\Dmap(V, W, \bfe, \bff) = (V, W, \bfe, \bff).
\]
Let $\Dmap(V, W, \bfe, \bff) = (\bfd, U)$, and let
$\Mmap(\bfd, U) = (V', W', \bfe', \bff')$.
Since
\[
\iota(V') = \iota(\delroot(t)) \cap U = \iota(V)
\]
and $\iota$ is injective,
it follows that $V' = V$.
The definition of $\bfd$ implies that exactly one of the
following possibilities occurs for $x\in X$.
\begin{eqnarray*}
\bfd(x) & = & \troot(e);\\
\bfd(x) & = &
\bfe(x)\in\nodes(e\hangon{1}\delroot(t)|V)\backslash\{\troot(e)\}
\mbox{ and }x\in W;\\
\bfd(x) & = & \bff(x)\in
\nodes(e\hangon{1}\delroot(t)|(\delroot(t)\backslash
V))\backslash\{\troot(e)\} \mbox{ and }x\notin W.
\end{eqnarray*}
Comparing this with the definition of $W'$ we see that
\begin{eqnarray*}
W' & = & \big\{\,x\in W \mid \bfe(x) = \troot(e)\,\big\} \cup\\
   &   & \qquad \big\{\,x\in W \mid \bfe(x)\in
\nodes(e\hangon{1}\delroot(t)|V)\backslash\{\troot(e)\}\,\big\}\\
   & = & W.
\end{eqnarray*}
It follows immediately from the definitions of $\bfd$, $\bfe$, and
$\bff$ that $\bfe'=\bfe$ and $\bff'=\bff$.
This completes the proof of the lemma.
\medbreak

\begin{lemma}\label{HopfCondLemma}
Let
$t\in\tree(\frX)$, $X\in\forest(\frX)$, $\bfd : X\rightarrow\nodes(t)$
be a function, $U\subset\delroot(t\hangon{d}X)$ be a subforest, and
\[
(V, W, \bfe, \bff) = \Mmap(\bfd, U).
\]
Then the underlying trees of
\[
e\hangon{1}\delroot(t\hangon{d}X)|U
\]
and
\[
(e\hangon{1}\delroot(t)|V)\hangon{e}X|W
\]
are equal.
\end{lemma}

\pf For simplicity we write $e\hangon{1}U$ for
$e\hangon{1}\delroot(t\hangon{d}X)|U$, etc.
We first observe that
\begin{eqnarray*}
\nodes(e\hangon{1}U) & = & \{\troot(e)\} \cup
\bigcup_{u\in U}\nodes(u)\\
   & = &\{\troot(e)\} \cup \bigcup_{u\in
U\cap\iota(\delroot(t))}\nodes(u) \cup {}\\
   &   &\qquad
 \bigcup_{u\in U\cap X}\nodes(u)\\
   & = &\{\troot(e)\} \cup \bigcup_{v\in V} \Biggl( \nodes(v)\cup
\bigcup_{\scriptstyle w\in X\atop\scriptstyle\bfd(w)\in\nodes(v)}
\nodes(w) \Biggr) \cup {}\\
   &   &\qquad \bigcup_{u\in U\cap X}\nodes(u)\\
   & = &\{\troot(e)\} \cup \bigcup_{v\in V}\nodes(v) \cup
\bigcup_{w\in W}\nodes(w)\\
   & = &\nodes\bigl((e\hangon{1}V)\hangon{e}W\bigr),
\end{eqnarray*}
so that the sets of nodes of the underlying trees are the same.

To complete the proof, we show that each node in these two trees has the
same parent in both trees.
If $\alpha$ is a node of $t$, then its parent in $t$ is the same as its
parent in $e\hangon{1}U$ and in $(e\hangon{1}V)\hangon{e}W$.
If $\alpha$ is a node of a tree $x\in X$, there are two possibilities.
In case $\alpha$ is not the root of $x$, then its parent in $x$ is the same
as its parent in $e\hangon{1}U$ and in $(e\hangon{1}V)\hangon{e}W$.
In case $\alpha$ is the root of $x$, its parent in $e\hangon{1}U$ is
$\bfd(x)$; its parent in $(e\hangon{1}V)\hangon{e}W$ is $\bfe(x)$.
Making the usual identification of $\troot(t)$ with $\troot(e)$, we have
that $\bfe(x)=\bfd(x)$, so that $\alpha$ has the same parent in both
trees in this case also.
This completes the proof of the lemma.
\medbreak

The next axiom says that the equality of trees in
Lemma~\ref{HopfCondLemma} holds in the family $\frX$.

\begin{axiom}\label{HopfCondAxiom}
Let $\frX$ be a family of trees, let
$t\in\tree(\frX)$, $X\in\forest(\frX)$, let
$\bfd : X\rightarrow\nodes(t)$
be a function, let $U\subset\delroot(t\hangon{d}X)$ be a subforest,
and let
\[
(V, W, \bfe, \bff) = \Mmap(\bfd, U).
\]
Then
\[
e\hangon{1}\delroot(t\hangon{d}X)|U
= (e\hangon{1}\delroot(t)|V)\hangon{e}X|W.
\]
\end{axiom}


\section{The Hopf algebra associated with a family}\label{HopfAlgSect}

In this section we describe the graded connected cocommutative Hopf
algebra associated with the family $\frX$,
and review some facts about the structure of
such Hopf algebras.

Suppose that $\frX$ is a family satisfying Axioms~1--9 of the
previous section.
Let $k$ be a field of characteristic~$0$.
Define $\kX$ to be the vector space over $k$ with basis $\tree(\frX)$.
We grade $\kX$ as follows:
if the underlying tree of $t\in\tree(\frX)$ has $n+1$ nodes, then $t$
has degree~$n$.
By Axiom~\ref{UnitAxiom}
there is only one $t\in\tree(\frX)$ whose underlying tree has one
node.
Therefore the graded vector space $\kX$ is connected.

We define a product on $\kX$ as follows:
if $t_1$, $t_2\in\tree(\frX)$, define
\refstepcounter{lemma}\label{ProdDef}
$$
t_1\cdot t_2 = \sum t_2\hangon{d}\delroot(t_1),
\eqno{\rm (\thelemma)}
$$
where the sum ranges over all possible linking maps
$\bfd : \delroot(t_1)\rightarrow\nodes(t_2)$.
We extend this product to all of $\kX$ by linearity.
It is immediate that this product respects the grading we have defined.
We now verify that this product is associative.
Note that
\[
(t_1\cdot t_2)\cdot t_3 =
\sum
t_3\hangon{e}\delroot(t_2\hangon{d}\delroot(t_1)),
\]
where the sum is taken over all pairs $(\bfd, \bfe)$ with
$\bfd : \delroot(t_1)\rightarrow\nodes(t_2)$, and
$\bfe : \delroot(t_2\hangon{d}\delroot(t_1))\rightarrow \nodes(t_3)$,
and that
\[
t_1\cdot(t_2\cdot t_3) =
\sum
(t_3\hangon{f}\delroot(t_2))\hangon{g}\delroot(t_1),
\]
where the sum is taken over all pairs ($\bff, \bfg)$ with
$\bff : \delroot(t_2)\rightarrow\nodes(t_3)$, and
$\bfg : \delroot(t_2)\rightarrow \nodes(t_3\hangon{f}\delroot(t_2))$.
Now Lemma~\ref{AssocPairsLemma} gives a one-one correspondence
between the set of pairs \{(\bfd, \bfe)\} over which the summation
equalling $(t_1\cdot t_2)\cdot t_3$ is taken,
and the set of pairs \{(\bff, \bfg)\} over which the summation
equalling $t_1\cdot (t_2\cdot t_3)$ is taken.
Axiom~\ref{EqualTreeAxiom} implies that the corresponding terms of the
summations are equal.
Axiom~\ref{HangonAxiom} implies that if $t\in\tree(\frX)$, and
$e\in\tree(\frX)$ is the unique element whose underlying tree has only
one node, then $t\cdot e = t$ and $e\cdot t = t$.

The definition of the product given in Equation~(\ref{ProdDef}) may
appear to be reversed.
The reason for this apparent reversal is that this product of trees is
one which has been used in applications involving data structures
representing differential operators (see \cite{Grossman}, \cite{FOCS},
\cite{NewDir}, and \cite{GLpara}).
The reversal is similar to the reversal which occurs in matrix
multiplication, in the correspondence between linear transformations on
a finite-dimensional vector space and matrices.

We now define a coproduct $\Delta : \kX\rightarrow\kX\otimes\kX$ as
follows:
if $t\in\tree(\frX)$ define
\[
\Delta(t) = \sum
(e\hangon{1}X|U)\otimes(e\hangon{1}X|(X\backslash U)),
\]
where $X=\delroot(t)$, and the sum is taken over all
subforests $U\subseteq X$.
($X\backslash Y$ denotes the set-theoretic difference of $X$ and $Y$.)
We extend $\Delta$ to all of $\kX$ by linearity.
It is immediate that this coproduct respects the grading on $\kX$.
We now verify that $\Delta$ is coassociative.
For trees, coassociativity is immediate:
all partitions of $\delroot(t)$ as a union of three disjoint
(possibly
empty) sets is achieved either by partitioning it into two disjoint
sets, and then partitioning the first set into two disjoint sets,
or
by partitioning it into two disjoint
sets, and then partitioning the second set into two disjoint sets.
Axiom~\ref{CoAssocAxiom} implies that this partitioning in two different
ways is equivalent in $\forest(\frX)$, and Axiom~\ref{HangonAxiom}
implies that we have
\[
(\Delta\otimes I)\circ\Delta = (I\otimes\Delta)\circ\Delta
\]
as maps from $\kX$ to $\kX\otimes\kX\otimes\kX$.
The counit $\epsilon : \kX\rightarrow k$ is defined as follows:
\[
\epsilon(t) = \left\{
\begin{array}{ll}
1 & \mbox{if $t = e$}\\
0 & \mbox{otherwise.}
\end{array}
\right.
\]
We extend $\epsilon$ to all of $\kX$ by linearity.
It follows that
\[
(\epsilon\otimes I)\circ\Delta = (I\otimes\epsilon)\circ\Delta = I
\]
from Axioms~\ref{TrivRestr} and~\ref{HangonAxiom}.
It follows from the fact that the set of all subsets of a set equals the
set of all complements of subsets of a set that this coalgebra is
cocommutative.

We now prove that the map $\Delta : \kX\rightarrow \kX\otimes\kX$ is an
algebra homomorphism, that is, that
\[
\Delta(t_1\cdot t_2) = \Delta(t_1)\cdot\Delta(t_2)
\]
for $t_1$, $t_2\in\tree(\frX)$.
We compute
\begin{eqnarray*}
\Delta(t_1\cdot t_2) & = &
\sum_{\bfd : X\rightarrow\nodes(t_2)}
\Delta(t_2\hangon{d}X)\\
& = & \sum
(e\hangon{1}Z_\bfd|U)\otimes(e\hangon{1}Z_\bfd|(Z_\bfd\backslash U)),
\end{eqnarray*}
where $Z_\bfd=\delroot(t_2\hangon{d}\delroot(t_1))$,
and the second sum ranges over all pairs $(\bfd, U)$, with
$\bfd : X\rightarrow\nodes(t_2)$ and
$U\subset \delroot(t_2\hangon{d}\delroot(t_1))$.
On the other hand
\begin{eqnarray*}
\Delta(t_1)\cdot\Delta(t_2) & = &
\left(\sum_{V\subset X}(e\hangon{1}X|V)\otimes
(e\hangon{1}X|(X\backslash V))\right)
\cdot\\
&   &\qquad
\left(\sum_{W\subset Y}(e\hangon{1}Y|W)\otimes
(e\hangon{1}Y|(Y\backslash W))\right)\\
& = & \sum
(e\hangon{1}Y|W)\hangon{e}X|V\otimes
(e\hangon{1}Y|(Y\backslash W))\hangon{f}X|(X\backslash V)),
\end{eqnarray*}
where $X=\delroot(t_1)$, $Y=\delroot(t_2)$, and
the last sum is taken over all quadruples $(V, W, \bfe, \bff)$ with
$V\subseteq X$, $W\subseteq Y$,
$\bfe : X|V\rightarrow\nodes(e\hangon{1}Y|W)$, and
$\bff : X|(X\backslash V)\rightarrow\nodes(e\hangon{1}Y|(Y\backslash W))$.
Now Lemma~\ref{HopfCondIndex} gives a one-one correspondence between the
terms of the summations equalling $\Delta(t_1\cdot t_2)$ and
$\Delta(t_1)\cdot\Delta(t_2)$.
Axiom~\ref{HopfCondAxiom} implies that the corresponding terms in the
summations are equal.

We summarize this discussion in the following theorem.
\begin{thm}
Let $\frX$ be a family of trees satisfying Axioms~1--9.
Then $\kX$ is a cocommutative graded connected Hopf algebra.
\end{thm}

If $A$ is a Hopf algebra, then the {\em primitive} elements of $A$ are
defined
\[
P(A) = \{\,a\in A \mid \Delta(a) = 1\otimes a + a \otimes 1\,\}.
\]
It can be shown that $P(A)$ is a Lie subalgebra of $A^{-}$,
which is the Lie algebra with the same underlying vector space as the
associative
algebra $A$, and in which the bracket operation is defined by
$[a,b] = ab - ba$.

If $L$ is a Lie algebra, then the universal enveloping algebra $U(L)$ is
a Hopf algebra.
If $x\in L$, then $\Delta(x) = 1\otimes x + x\otimes 1$, and
$\epsilon(x) = 0$.
The maps $\Delta$ and $\epsilon$ are extended to all of $U(L)$ using the
facts that
$\Delta$ is an algebra homomorphism, and that $L$ generates $U(L)$ as an
algebra.
The following theorem gives a basis for $U(L)$ in terms of an ordered
basis for $L$.

\begin{thm}[Poincar\'e-Birkhoff-Witt]\label{PBWThm}
Let $L$ be a Lie algebra with ordered basis
$x_1$, \ldots, $x_n$, \ldots.
Then
\[
\{\,x^{e_1}_{i_1}\cdots x^{e_t}_{i_t}
\mid i_1 < \cdots < i_t; \;\; e_k > 0\,\}
\]
is a basis for $U(L)$.
\end{thm}

See \cite[page~159]{Jac:LAs} for a proof.

\begin{thm}[Milnor-Moore]\label{MMThm}
Let $A$ be a cocommutative graded connected Hopf algebra.
Then
\[
A \cong U(P(A))
\]
as Hopf algebras.
\end{thm}

See \cite[page~244]{MM} or \cite[page~274]{Moss} for a proof.
\medbreak

If $X$ is a set, denote by $k{<}X{>}$ the free associative algebra over
$k$ generated by $X$.
Then $k{<}X{>}$ is a cocommutative Hopf algebra, with
$\Delta(x) = 1\otimes x + x\otimes 1$ for $x\in X$.
It can be shown that $P(k{<}X{>})$ is the free Lie algebra generated by
$X$.


\section{The family of labeled trees}\label{LTSect}

In this section we discuss the structure of
$k\{\frLT(E_1$, \ldots, $E_M)\}$.
Note that if $M = 1$ we are essentially discussing $k\{\frT\}$.
We give a description of $P(k\{\frLT(E_1$, \ldots, $E_M)\})$
and use it, together with the Milnor-Moore Theorem and the
Poincar\'e-Birkhoff-Witt Theorem, to give a new proof of the recurrence
relation
for the number of rooted trees with $n$ nodes first given by Cayley
\cite{Cayley} in~1857.

Let $\frLT_1(E_1$, \ldots, $E_M)$ be the set of labeled trees
$t\in\frLT(E_1$, \ldots, $E_M)$ whose root has only one child.

\begin{thm}\label{prims}
The set $\frLT_1(E_1$, \ldots, $E_M)$ is a basis for
$P(k\{\frLT(E$, \ldots, $E_M)\})$.
\end{thm}

\pf Denote $\frLT(E_1$, \ldots, $E_M)$ by $\frLT$,
and $\frLT_1(E_1$, \ldots, $E_M)$ by $\frLT_1$.
It is easily checked that if $t\in\frLT_1$,
then $t\in P(k\{\frLT\})$.
We now show that the elements of $\frLT_1$ span $P(k\{\frLT\})$.
Define
\[
\pi : k\{\frLT\}\otimes k\{\frLT\} \rightarrow k\{\frLT\}
\]
as follows:
if $t_1$, $t_2\in\frLT$, let $\pi(t_1\otimes t_2)$ be the element of
$\frLT$ formed by identifying the roots of $t_1$ and $t_2$.
It is easily checked that if the root of $t\in\frLT$ has $r$ children,
then $\pi\circ\Delta(t) = 2^rt$.
On the other hand,
if $a = \sum a_t t\in P(k\{\frLT\})$, then $\pi\circ\Delta(a) = 2a$.
Since the elements of $\frLT$ are linearly independent, it follows that
$a_t=0$ if the root of $t$ has more than one child.
This completes the proof of the theorem.
\medbreak

\begin{thm}\label{CayleyTh}
Let $t_n$ be the number of rooted trees with $n$ nodes,
with all nodes but the root labeled using the formal symbols
$\{E_1$, \ldots, $E_M\}$.
Then
\begin{eqnarray*}
t_1 & = & 1\\
t_{n+1} & = & \sum_{m_1 + 2m_2 +\cdots+rm_r = n}
{m_1 + Mt_1 - 1\choose m_1}\cdots{m_r + Mt_r - 1\choose m_r}.
\end{eqnarray*}
\end{thm}

\pf Let
\[
a_n = \dim {k\{\frLT(E_1, \ldots, E_M)\}}_n
\]
and let
\[
p_n = \dim {P(k\{\frLT(E_1, \ldots, E_M)\})}_n.
\]
The definition of the grading on $\kX$ implies that $t_{n+1} = a_n$.
Theorem~\ref{prims} implies that $p_n = Ma_{n-1}$.
Since the number of monomials of length $m$ of the form
$x_1^{e_1}\cdots x_{p_t}^{e_{p_t}}$,
where $\{x_i\}$ is an ordered basis of
${P(k\{\frLT(E_1, \ldots, E_M)\})}_t$, is
$m + p_t - 1\choose m$, it follows that
\begin{eqnarray*}
a_{n} & = & \sum_{m_1 + 2m_2 +\cdots+rm_r = n}
{m_1 + p_1 - 1\choose m_1}\cdots{m_r + p_r - 1\choose m_r}\\
      & = & \sum_{m_1 + 2m_2 +\cdots+rm_r = n}
{m_1 + Ma_0 - 1\choose m_1}\cdots{m_r + Ma_{r-1} - 1\choose m_r}.
\end{eqnarray*}
The statement of the theorem follows immediately from this.
\medbreak

The following result (with $M=1$) was proved by Cayley~\cite{Cayley} in
1857.

\begin{corollary}
Let $T(z) = \sum_{n = 1}^\infty t_n z^n$.
Then
\[
T(z) = z \prod_{n=1}^\infty {(1 - z^n)}^{-Mt_n}.
\]
\end{corollary}

\pf This follows immediately from the Theorem~\ref{CayleyTh}
upon observing that
\[
{r + k - 1\choose k} = (-1)^k{-r\choose k}.
\]


\section{The family of labeled ordered trees}\label{LOTSect}

In this section we discuss the structure of
$k\{\frLOT(E_1$, \ldots, $E_M)\}$.
We will show that this Hopf algebra is isomorphic to the free
associative algebra generated by
$\frLOT_1(E_1$, \ldots, $E_M)$, the set of labeled ordered rooted trees
whose root has exactly one child.
This fact allows us to give a recurrence relation for the number of
labeled ordered rooted trees.
This recurrence can be solved to get the number of ordered rooted trees
with $n$ nodes.
This number can be shown~\cite{KnuthI} to be the same as the number of
binary trees with $n-1$ nodes, which was given by Catalan~\cite{Catalan}
in~1838, and by Cayley~\cite{Cayley2} in~1859.

\begin{thm}\label{freeOT}
\[
k\{\frLOT(E_1, \ldots, E_M)\} \cong k{<}\frLOT_1(E_1, \ldots, E_M){>}.
\]
\end{thm}

\pf Write $\frLOT$ for $\frLOT(E_1$, \ldots, $E_M)$,
and $\frLOT_1$ for $\frLOT_1(E_1$, \ldots, $E_M)$
We introduce a filtration on $\kFLOT$ by defining $F_p\kFLOT$ to be the
subspace of $\kFLOT$ spanned by all monomials of length $\le p$.
It is clear that
\begin{eqnarray*}
F_{-1}\kFLOT & = & 0,\\
\bigcup_pF_p\kFLOT & = & \kFLOT,\\
(F_p\kFLOT)\cdot(F_q\kFLOT) & \subseteq & F_{p+q}\kFLOT.
\end{eqnarray*}
We introduce a filtration on $\kLOT$ by defining $F_p\kLOT$ to be the
subspace spanned by all trees whose root has $p$ or fewer children.
This filtration satisfies
\begin{eqnarray*}
F_{-1}\kLOT & = & 0,\\
\bigcup_pF_p\kLOT & = & \kLOT,\\
(F_p\kLOT)\cdot(F_q\kLOT) & \subseteq & F_{p+q}\kLOT.
\end{eqnarray*}
More precisely, if the root of $t_1$ has $p$ children, and the root of
$t_2$ has $q$ children, then
\[
t_1\cdot t_2 = \sum_{\bfd : \delroot(t_1)\rightarrow\nodes(t_2)}
t_2\hangon{d}\delroot(t_1),
\]
and the root of $t_2\hangon{d}\delroot(t_1)$ has $q+r$ children, where
$r$ is the number of $x\in\delroot(t_1)$ satisfying
$\bfd(x) = \troot(t_2)$.
Let ${\bf d_0} : \delroot(t_1)\rightarrow\nodes(t_2)$ be defined by
${\bf d_0}(x) = \troot(t_2)$ for all $x\in\delroot(t_1)$.
Then
\[
t_1\cdot t_2 = t_2\hangon{d_0}\delroot(t_1) + \tau,
\]
where $\tau\in F_{p+q-1}\kLOT$.

There is a unique algebra homomorphism $\phi : \kFLOT\rightarrow\kLOT$
which is the identity on $\frLOT_1\subset\frLOT$.
It is clear that $\phi(F_p\kFLOT)\subseteq F_p\kLOT$.
More precisely, if $x_1$, \ldots, $x_p\in\frLOT_1$,
then $\phi(x_1\cdots x_p) = \phi(x_1)\cdots\phi(x_p)$ is congruent
modulo $F_{p-1}\kLOT$ to the tree whose root has exactly the following
$p$ children:
the children of the roots of $x_p$, \ldots, $x_1$ in that order.
Therefore $\phi$ induces an isomorphism
\[
F_p\kFLOT/F_{p-1}\kFLOT\rightarrow F_p\kLOT/F_{p-1}\kLOT.
\]
It follows that $\phi$ is an isomorphism.
This completes the proof of the theorem.
\medbreak

\begin{thm}\label{OrdTree}
Let $t_n$ be the number of ordered rooted trees with $n$ nodes,
with all of the nodes but the root labeled using the formal symbols
$\{E_1$, \ldots, $E_M\}$.
Then
\begin{eqnarray*}
t_1     & = & 1\\
t_{n+1} & = & M(t_1t_n + t_2t_{n-1} + \cdots + t_nt_1).
\end{eqnarray*}
\end{thm}

\pf Let
\[
a_n = \dim {k\{\frLOT(E_1, \ldots, E_M)\}}_n,
\]
and
\[
g_n = \card \{\,t\in\frLOT_1(E_1, \ldots, E_M) \mid
      \mbox{$t$ has $n+1$ nodes}\,\},
\]
the number of free generators of degree $n$.
The definition of the grading on $\kX$ implies that $t_{n+1}=a_n$.
It follows from the fact that $\kLOT\cong\kFLOT$ that $g_n = Ma_{n-1}$,
so
\begin{eqnarray*}
a_n & = & g_1a_{n-1} + \cdots + g_na_0\\
    & = & M(a_0a_{n-1} + \cdots + a_{n-1}a_0),
\end{eqnarray*}
or
\[
t_{n+1} = M(t_1t_n + \cdots + t_nt_1).
\]
This completes the proof of the theorem.
\medbreak

\begin{corollary}
Let $t_n$ be the number of ordered rooted trees with $n$ nodes,
with all of the nodes but the root labeled using the formal symbols
$\{E_1$, \ldots, $E_M\}$.
Then
\[
t_n = {M^{n-1}\over n}{2n-2 \choose n-1}.
\]
\end{corollary}

\pf The solution of the recurrence in Theorem~\ref{OrdTree} is adapted
from \cite[pages 388--389]{KnuthI}.
This proof first appears in~\cite{Binet}.
Let
\[
T(z) = \sum_{n=1}^\infty t_nz^n.
\]
Then Theorem~\ref{OrdTree} implies that
\[
T(z) = M{T(z)}^2 + z
\]
which, together with the fact that the coefficient of $z$ in
$T(z)$ equals $1$ implies that
\[
T(z) = {1 - \sqrt{1-4Mz}\over 2M}.
\]
Now
\[
\sqrt{1-4Mz} = \sum_{n=0}^\infty{(-1)}^n4^nM^n{{1\over2}\choose n}z^n
\]
and
\[
{{1\over2}\choose n} = {{(-1)}^{n-1}\over2^{2n}(2n-1)}{2n\choose n},
\]
so
\[
T(z) = \sum_{n=1}^\infty {M^{n-1}\over n}{2n-2\choose n-1}z^n,
\]
which proves the corollary.


\section{Heap-ordered trees and other families}\label{XSect}

In this section, we describe some other families of trees related to
$\frOT$, and prove a generalization of Theorem~\ref{freeOT} for these
families.
We then apply this result to get a recurrence relation on the number of
indecomposable heap-ordered trees.
First we give a description of the category of families of trees.

If $\frX$ and $\frak Y$ are families of trees, we define a morphism
$f : \frX\rightarrow {\frak Y}$ to be a pair of maps
$\tree(\frX)\rightarrow\tree({\frak Y})$ and
$\forest(\frX)\rightarrow\forest({\frak Y})$
which are the identity on the underlying trees and forests,
and which commute with $\delroot$, $|$, and $\leftharpoonup$.
We thus define the category $\bf Family$ of families of trees.
Note that $\frT$ is a terminal object in $\bf Family$.
In Section~\ref{HopfAlgSect}
we defined the functor $k\{{-}\}$ from $\bf Family$ to the category of
cocommutative graded connected Hopf algebras over $k$.

In the example below, we will speak of generations of a node in a tree:
the singleton $\{\alpha\}$ is the first generation of $\alpha$;
if $\{\beta_1$, \ldots, $\beta_k\}$ is the $n^{\rm th}$ generation of
the node $\alpha$, then the ${n+1}^{\rm st}$ generation of $\alpha$ is
the set of children of the nodes $\beta_1$, \ldots, $\beta_k$.

\begin{example}
The family $\frGT_n$ of $n$-generation ordered trees.
The set $\tree(\frGT_n)$ is the set of finite rooted trees,
together with, for each node $\alpha$, a total ordering on the set
$X_\alpha$ consisting of the nodes which constitute the first $n$
generations of $\alpha$.
The ordering on $X_\alpha\cap X_\beta$ must be the same, whether induced
as a subset of $X_\alpha$ or of $X_\beta$.
The orderings must also satisfy the condition that each node $\alpha$
precede all of its descendants in $X_\alpha$.

The set $\forest(\frGT_n)$ is the set of finite forests of finite rooted
trees,
together with, for each node $\alpha$, a total ordering on the set
$X_\alpha$ consisting of the nodes which constitute the first $n$
generations of $\alpha$,
and a total ordering on the set $Y$ consisting of the nodes which
constitute the first $n-1$
generations of the nodes which
are the roots of the trees in the forest.
The ordering on $X_\alpha\cap X_\beta$  (or on $X_\alpha\cap Y$) must be
the same, whether induced as a subset of $X_\alpha$ or of $X_\beta$ (or
of $Y$).
The orderings must also satisfy the condition that each node $\alpha$
precede all of its descendants in $X_\alpha$.

The map $\delroot$ sends each tree into the forest formed by deleting
its root, with the orderings on the sets of nodes given in the obvious
way.

If $V\in\forest(\frGT_n)$ and $U\subseteq V$ is a subforest, then $V|U$
is the forest $U$, with the same orderings on the sets $X_\alpha$ as in
$V$, and with the ordering on the set $Y_U$ associated with $U$ being
the one induced by the fact that it is a subset of the set $Y_V$
associated with $V$.

If $t\in\tree(\frGT_n)$, $U\in\forest(\frGT_n)$, and
$\bfd : U\rightarrow\nodes(t)$,
then $t\hangon{d}U$ is the tree formed by linking the root of each tree
in $U$ to the node $\bfd(u)$ of $t$.
If $\alpha$ was originally a node of $t$, the new set $X_\alpha$ is
ordered as follows:
the original descendants of $\alpha$ in $t$, up to the $n^{\rm th}$
generation, preserve their original order in the new $X_\alpha$;
the descendants of $\alpha$, up to the $n^{\rm th}$ generation,
which were among the nodes of trees in $U$, preserve the
order they originally had in $Y_U$;
all of the original nodes of $t$ which are descendants of $\alpha$,
up to the $n^{\rm th}$ generation, precede all of the descendants of
$\alpha$, up to the $n^{\rm th}$ generation, which were originally nodes
of trees in $U$.
If $\alpha$ was originally a node of a tree in $U$, the set $X_\alpha$
is unchanged.
\end{example}

\begin{example}
The family $\frLGT_n$ of labeled $n$-generation ordered trees.
This is formed by extending the definition of the family $\frGT_n$
in the same way that the definition of $\frT$ was extended in
Example~\ref{LTEx} to a definition of $\frLT$.
\end{example}

In some cases, the family $\frGT_n$ is isomorphic to a previously
defined family.
The family $\frGT_1$ is isomorphic to the family $\frT$.
The family $\frGT_2$ is isomorphic to the family $\frOT$
(the map $\frOT\rightarrow\frGT_2$ is given by:
if $t\in\frOT$, construct orderings on the sets of nodes and their
children by requiring that each node precede all of its children, and
keep the ordering of the children unchanged; the inverse map
$\frGT_n\rightarrow\frOT$ is given by:
if $t\in\frGT_2$, order the children of each node in the same way that
they are ordered in $\frGT_2$).

There is a morphism $\frGT_n\rightarrow\frGT_{n-1}$ given by restricting
the ordering on the set of the first $n$ generations of each node
$\alpha$ to the first $n-1$ generations.
These morphisms are surjections.
It is clear that $\frHOT = \lim\limits_{\longleftarrow}\frGT_n$.

Let $\frX$ be a family of trees.
If $t_1$, $t_2\in\tree(\frX)$, write $t_1\odot t_2$ for the tree
\[
t_2\hangon{d_0}\delroot(t_1),
\]
where ${\bf d_0} : \delroot(t_1)\rightarrow\nodes(t_2)$ is the function
defined by ${\bf d_0}(w)=\troot(t_2)$, for all $w\in\delroot(t_1)$.
Note that Axiom~\ref{EqualTreeAxiom} implies that $\odot$ is
associative.
If $t\in\tree(\frX)$ with $t\not=e$, we say that $t$ is
{\em indecomposable} if $t=t_1\odot t_2$ implies that $t_1=e$ or
$t_2=e$.
Let
\[
I(\frX) = \{\,t\in\tree(\frX) \mid
\mbox{$t$ is indecomposable}\,\}.
\]
By finiteness, every tree $t\in\tree(\frX)$ can be written
\[
t = t_1\odot\cdots\odot t_p,
\]
with $t_i\in I(\frX)$.

For example, if $\frX=\frT$ or $\frOT$, the indecomposable trees are
those trees whose roots have exactly one child.
In $\frHOT$, writing $t=t_1\odot t_2$ corresponds to partitioning
$\delroot(t)$ into two disjoint subforests $V=\delroot(t_1)$ and
$W=\delroot(t_2)$,
such that every node of every
tree in $V$ precedes every node of every tree in $W$.
An indecomposable tree $t\in\tree(\frHOT)$ is therefore one for which
$\delroot(t)$ cannot be so partitioned.
There is an analogous (but not so clearly expressible) description of
indecomposable trees in $\frGT_n$ for $n>2$.

The family $\frX$ is said to have the {\em unique decomposition
property} if, whenever
\[
x_1\odot\cdots\odot x_p = y_1\odot\cdots\odot y_q,
\]
with $x_i$, $y_j\in I(\frX)$, the ordered sequences $(x_1$, \ldots,
$x_p)$ and $(y_1$, \ldots, $y_q)$ are equal.
Essentially, the unique decomposition property is a non-commutative
unique factorization condition.
Note that the families $\frOT$, $\frHOT$, and $\frGT_n$, for $n\ge2$,
all have the unique decomposition property.
The family $\frT$ does not have the unique decomposition property,
because $\odot$ is commutative in $\frT$.

The following theorem can be thought of as a generalization of
Theorem~\ref{freeOT}.

\begin{thm}\label{freeUDP}
Let $\frX$ be a family of trees with the unique decomposition property.
Then
\[
\kX \cong \kIX.
\]
\end{thm}

\pf We introduce a filtration on $\kIX$ by filtering $I(\frX)$:
if $t\in I(\frX)$, we say that $t\in F_pI(\frX)$ if the root of $t$ has
$\le p$ children.
Next define $F_p\kIX$ to be the subspace spanned by the monomials
$x_1\cdots x_k$, where $x_i\in F_{p_i}I(\frX)$ and
$p_1+\cdots+p_k \le p$.
It is clear that
\begin{eqnarray*}
F_{-1}\kIX & = & 0,\\
\bigcup_pF_p\kIX & = & \kIX,\\
(F_p\kIX)\cdot(F_q\kIX) & \subseteq & F_{p+q}\kIX.
\end{eqnarray*}

We next filter $\kX$ by defining $F_p\kX$ to be the subspace spanned
by the $t\in\frX$ whose roots have $\le p$ children.
This filtration satisfies
\begin{eqnarray*}
F_{-1}\kX & = & 0,\\
\bigcup_pF_p\kX & = & \kX,\\
(F_p\kX)\cdot(F_q\kX) & \subseteq & F_{p+q}\kX.
\end{eqnarray*}

There is a unique algebra homomorphism
$\phi : \kIX\rightarrow\kX$ which is the identity on
$I(\frX)\subset\frX$.
It is easily checked that $\phi(F_p\kIX)\subseteq F_p\kX$.
Therefore $\phi$ induces a map
\[
\overline\phi : F_p\kIX/F_{p-1}\kIX\rightarrow F_p\kX/F_{p-1}\kX.
\]
If $x_1$, \ldots, $x_k\in I(\frX)$, $x_1\cdots x_k\in F_p\kIX$, and
$x_1\cdots x_k\notin F_{p-1}\kIX$, then $\phi(x_1\cdots x_k) =
\phi(x_1)\cdots\phi(x_k)$ is congruent modulo $F_{p-1}\kX$ to
$x_1\odot\cdots\odot x_k$.
Since $\frX$ has the unique decomposition property, the elements
$x_1\odot\cdots\odot x_k$ for which $x_1\cdots x_k\in F_p\kIX$, and
$x_1\cdots x_k\notin F_{p-1}\kIX$, are linearly independent modulo
$F_{p-1}\kX$.
Therefore, the map $\overline\phi$ is injective.
It follows from the fact that every $x\in\frX$ can be written
$x = x_1\odot\cdots\odot x_k$,
with $x_i\in I(\frX)$, that $\overline\phi$ is surjective.
It follows that $\phi$ is an isomorphism.
This completes the proof of the theorem.
\medbreak

It can be shown that a Hopf algebra structure can also be defined on
$\kX$ using $\odot$ as product, rather than the product defined in
Section~\ref{HopfAlgSect}.
This Hopf algebra is just the associated graded Hopf algebra which
is constructed from the filtration $F_p\kX$ defined in the proof of
Theorem~\ref{freeUDP}.
We summarize this in the following corollary.
\begin{corollary}
The vector space $\kX$, with the coproduct and counit defined in
Section~\ref{HopfAlgSect}, and with product $\odot$, is a cocommutative
bigraded connected Hopf algebra which is isomorphic to
\[
\sum_{p\ge0}F_p\kX/F_{p-1}\kX.
\]
\end{corollary}
\medbreak

We now use Theorem~\ref{freeUDP} to give a recurrence relation on the
number of indecomposable heap-ordered trees.
We first compute the number of heap-ordered trees.
\begin{lemma}\label{HOTcnt}
The number of distinct heap-ordered trees with $n+1$ nodes is $n!$.
\end{lemma}

\pf The assertion of the lemma is clear for $n=0$.
If $n>0$, a heap-ordered tree with $n+1$ nodes can be formed by linking
the ${n+1}^{\rm st}$ node, which is last in the total ordering on the
nodes, to any of the $n$ nodes in any of the $(n-1)!$ heap-ordered trees
with $n$ nodes.
It is clear that the resulting $n!$ heap-ordered trees are all distinct.
This completes the proof of the lemma.
\medbreak

\begin{corollary}
The number $g_n$ of distinct indecomposable heap-ordered trees with
$n+1$ nodes satisfies
\begin{eqnarray*}
g_1 & = & 1\\
g_n & = & n! - (n-1)!\,g_1 - \cdots - 1!\,g_{n-1}.
\end{eqnarray*}
\end{corollary}

\pf Let $t_n$ be the number of heap-ordered trees with $n+1$ nodes.
By Theorem~\ref{freeUDP}, the algebra $\kHOT$ is freely generated by
$I(\frHOT)$, so
\[
t_n = g_1t_{n-1} + \cdots + g_nt_0.
\]
Since $t_n = n!$ by Lemma~\ref{HOTcnt}, the corollary follows.


\begin{thebibliography}{99}
\bibitem{Binet} M. J. Binet, {\em R\'eflexions sur le probl\`eme de
d\'eterminer le nombre de mani\`eres dont une figure rectiligne peut
\^etre partag\'ee en triangles au moyen de ses diagonals},
J. Math. Pures Appl. (1) {\bf4} (1839), 79--90.
\bibitem{Catalan} E. Catalan, {\em Note sur une \'equation aux
diff\'erences finies},
J. Math. Pures Appl. (1) {\bf3} (1838), 508--516.
\bibitem{Cayley} A. Cayley, {\em On the theory of the analytical forms
called trees}, in ``Collected Mathematical Papers of Arthur
Cayley,'' Cambridge Univ. Press, Cambridge, 1890, {\bf 3}, 242--246.
\bibitem{Cayley2} A. Cayley, {\em On the analytical forms called trees.
Second part}, in ``Collected Mathematical Papers of Arthur
Cayley,'' Cambridge Univ. Press, Cambridge, 1891, {\bf 4}, 112--115.
\bibitem{Grossman} R. Grossman, {\em Evaluation of expressions involving
higher order derivations,}
Center for Pure and Applied Mathematics,
PAM--367, University of California, Berkeley.
\bibitem{FOCS} R. Grossman and R. G. Larson,
{\em Labeled trees and the algebra of differential operators,}
Center for Pure and Applied Mathematics,
PAM--368, University of California, Berkeley.
\bibitem{NewDir} R. Grossman and R. G. Larson,
{\em The symbolic computation of higher-order derivations:
symmetries of expressions and actions of group algebras,}
Center for Pure and Applied Mathematics,
PAM--370, University of California, Berkeley.
\bibitem{GLpara} R. Grossman and R. G. Larson,
{\em Labeled trees and the efficient computation of
derivations}, in preparation.
\bibitem{Jac:LAs} N. Jacobson, ``Lie algebras,'' Interscience,
New York, 1962.
\bibitem{JoniRota} S. A. Joni and G.-C. Rota,
{\em Coalgebras and bialgebras in combinatorics},
Stud. Appl. Math. {\bf61} (1979), 93--139; reprinted
in ``Umbral Calculus and Hopf Algebras,''
Amer. Math. Soc., Providence, 1982,
1--47.
\bibitem{KnuthI} D. E. Knuth, ``The Art of Computer Programming, Vol. I,
Fundamental Algorithms,'' $2^{\rm nd}$ Ed., Addison-Wesley, Reading,
1973.
\bibitem{MM} J. W. Milnor and J. C. Moore, {\em On the structure of Hopf
algebras,} Ann. Math. (2) {\bf 81} (1965), 211--264.
\bibitem{NiSw} W. Nichols and M. Sweedler,
{\em Hopf algebras and combinatorics},
in ``Umbral Calculus and Hopf Algebras,''
Amer. Math. Soc., Providence, 1982,
49--84.
\bibitem{Bodo} B. Pareigis, ``Four Lectures on Hopf Algebras,''
Centre de Recerca Matem\`atica, Institut d'Estudis Catalans,
No.~6, Octubre 1984.
\bibitem{Moss} M. Sweedler, ``Hopf Algebras,'' W. A. Benjamin, New York,
1969.
\bibitem{Tarjan} R. E. Tarjan, ``Data Structures and Network
Algorithms,'' SIAM, Philadelphia, 1983.
\end{thebibliography}
\end{document}